\documentclass[reqno]{article}
\usepackage{amsmath, amsfonts, amscd, eucal}
\usepackage{latexsym}
\usepackage{amssymb}
\newcommand{\be}{\begin{equation}}
      \newcommand{\ee}{\end{equation}}
      \newcommand{\ba}{\begin{eqnarray}}
       \newcommand{\ea}{\end{eqnarray}}
\newcommand{\ban}{\begin{eqnarray*}}
\newcommand{\ean}{\end{eqnarray*}}

\newcommand{\pt}{\partial}
\newcommand{\lp}{\langle}
\newcommand{\rp}{\rangle}
\newcommand{\ra}{\rightarrow}

 \newcommand{\qed}{\hspace*{\fill}\rule{3mm}{3mm}\quad \vspace{.2cm}}
 \newcommand{\Pf}{\noindent {\bf Proof:} }
 \newcommand{\Rk}{\noindent {\bf Remark} }

\newcommand{\sect}[1]{\section{#1} \setcounter{equation}{0}}

\newtheorem{theo}{Theorem}[section]

\begin{document}
\newtheorem{defn}[theo]{Definition}
\newtheorem{ques}[theo]{Question}
\newtheorem{lem}[theo]{Lemma}
\newtheorem{prop}[theo]{Proposition}
\newtheorem{coro}[theo]{Corollary}
\newtheorem{ex}[theo]{Example}
\newtheorem{note}[theo]{Note}

\title{On the Stability of Riemannian Manifold with Parallel Spinors}
\author{Xianzhe Dai\thanks {Math Dept, UCSB, Santa Barbara, CA 93106 \tt{Email:
dai@math.ucsb.edu}.} \and Xiaodong
Wang\thanks{Department of Mathematics, MIT, Cambridge, MA 02139,
USA. \tt{Email:xwang@math.mit.edu}} \and Guofang Wei\thanks {Math
Dept. UCSB. \tt{Email: wei@math.ucsb.edu}. Partially supported by
NSF Grant \# DMS-0204187.}}
\date{\it Dedicated to Jeff Cheeger for his sixtieth birthday} \maketitle

\begin{abstract}Inspired by the recent work \cite{HHM}, we prove two
stability results for compact Riemannian manifolds with nonzero
parallel spinors. Our first result says that Ricci flat metrics
which also admits nonzero parallel spinors are stable (in the
direction of changes in conformal structures) as the critical
points of the total scalar curvature functional. In fact, we show
that the Lichnerowicz Laplacian, which governs the second
variation, is the square of a twisted Dirac operator. Our second
result, which is a local version of the first one, shows that any
metrics of positive scalar curvature cannot lie too close to a
metric with nonzero parallel spinor. We also prove a rigidity
result for special holonomy metrics. In the case of $SU(m)$
holonomy, the rigidity result implies that scalar flat
deformations of Calabi-Yau metric must be Calabi-Yau. Finally we
explore the connection with a positive mass theorem of \cite{d},
which presents another approach to proving these stability and
rigidity results.
\end{abstract}

\newcommand{\inj}{\mbox{inj}}
\newcommand{\vol}{\mbox{vol}}
\newcommand{\diam}{\mbox{diam}}
\newcommand{\Ric}{\mbox{Ric}}
\newcommand{\Iso}{\mbox{Iso}}
\newcommand{\Hess}{\mbox{Hess}}
\newcommand{\divg}{\mbox{div}}
\newcommand\grd{\nabla}
\newcommand{\Dirac}{\mathcal{D}}
\newcommand{\M}{\mathcal{M}}
\newcommand{\lx}{L_X\raisebox{0.5ex}{$g$}}
\newcommand{\cuv}{\overset{\hspace{.5ex}\circ}{R}}
\newcommand{\Sr}{\mathcal{S}}
\newcommand{\ts}{\otimes}
\newcommand{\Mh}{\hat{M}}
\newcommand{\gh}{\hat{g}}

\def\operatorname#1{{\rm #1\,}}
\def\op{\operatorname}
\def\lam{\lambda}
\def\s{\sigma}
\def\ph{\phi}
\def\ps{\psi}
\def\e{\epsilon}
\def\d{\delta}
\def\D{\Delta}

\sect{Introduction}

One of the most fruitful approaches to finding the `best' (or
canonical) metric on a manifold has been through the critical
points of a natural geometric functional. In this approach one is
led to the study of variational problems and it is important to
understand the stability issue associated to the variational
problem. Consider the space $\M$ of Riemannian metrics on a
compact manifold $M$ (our manifolds here are assumed to have empty
boundary). It is well known that the critical points of the total
scalar curvature functional (also known as the Hilbert-Einstein
action in general relativity) are Ricci flat metrics. It is also
well known that the total scalar curvature functional behaves in
opposite ways along the conformal deformations and its transversal
directions (i.e., when the conformal structure changes). The
variational problem in the conformal class of a metric (volume
normalized) is the famous Yamabe problem, which was resolved by Aubin
and Schoen. In this paper we study the
stability for the total scalar curvature functional when we
restrict to the transversal directions, that is, the space of
conformal structures.

This has to do with the second variation of the total scalar
curvature functional restricted to the traceless transverse
symmetric $2$-tensors, which is given in terms of Lichnerowicz
Laplacian \cite{Besse}. Our first result shows that Riemannian
manifolds with nonzero paprallel spinors (which are necessarily
Ricci flat) are stable in this sense. In fact, by identifying the
symmetric $2$-tensors with twisted spinors, we show that the
Lichnerowicz Laplacian is the square of a twisted Dirac operator,
hence positive semi-definite\footnote{It has since been brought to
our attention that this is already implicit in the work of
McKenzie Wang \cite{Wa2}, where the general case of Killing spinor
is discussed. We thank McKenzie for bringing this work to our
attention.}.

\begin{theo} If a compact Riemannian manifold $(M, g)$ has a cover which is spin
and admits nonzero parallel spinors, then the Lichnerowicz
Laplacian $\mathcal{L}_g$ is positive semi-definite.
\end{theo}

This settles an open question raised in \cite{KW} about thirty
years ago in the case when the Ricci flat manifold has a spin
cover with nonzero parallel spinors. It is an interesting open
question of how special our metric is compared to the general
Ricci flat metric (Cf. Section 5) but we note that, so far, all
known examples of compact Ricci flat manifolds are of this type,
namely, they admit a spin cover with nonzero parallel spinors.

We then prove that, in fact, there exists a neighborhood of the
metric with nonzero parallel spinors, which contains no metrics
with positive scalar curvature. This can be thought of as a local
version of our previous (infinitesimal) stability result.

\begin{theo} Let $(M, g)$ be a compact Riemannian manifold
which admits a spin cover with nonzero parallel spinors. Then $g$
cannot be deformed to positive scalar curvature metrics. By this
we mean that there exists no path of metrics $g_t$ such that
$g_0=g$ and the scalar curvature $S(g_t)>0$ for $t>0$. In fact, if
$(M,g)$ is simply connected and irreducible, then there is a
neighborhood of $g$ in the space of metrics which does not contain
any metrics with positive scalar curvature.
\end{theo}

The existence of metrics with positive scalar curvature is a well
studied subject, with important work such as \cite{L}, \cite{H},
\cite{SY}, \cite{GL}, culminating in the solution of the
Gromov-Lawson conjecture in \cite{S} for simply connected
manifolds. Thus, a $K3$ surface does not admit any metric of
positive scalar curvature, but a simply connected Calabi-Yau
$3$-fold does. Of course, a Calabi-Yau admits nonzero parallel
spinors (in fact, having nonzero parallel spinors is more or less
equivalent to having special holonomy except quaternionic
K\"ahler, cf. Section 3). Thus on a Calabi-Yau $3$-fold there
exists metrics of positive scalar curvature but they cannot be too
close to the Calabi-Yau metric.

One should also contrast our result with an old result of
Bourguignon, which says that an metric with zero scalar curvature
but nonzero Ricci curvature can always be deformed to a metric
with positive scalar curvature (essentially by Ricci flow).

The local stability theorem implies a rigidity result, Theorem
3.4, from which we deduce the following interesting application.

\begin{theo} Any scalar flat deformation of a Calabi-Yau metric on
a compact manifold must be Calabi-Yau. The same is true for the
other special holonomy metrics, i.e., hyperk\"ahler, $G_2$, and
$Spin(7)$.
\end{theo}

Our result generalized a theorem in \cite{Wa2} for Einstein
deformations.

The proof of our second result, the local stability theorem,
actually gives a very nice picture of what happens to scalar
curvature near a metric with parallel spinor (we'll call it
special holonomy metric): one has the finite dimensional smooth
moduli of special holonomy metrics; along the normal directions,
the scalar curvature of the Yamabe metric in the conformal class
will go negative. In other words, we have

\begin{theo} The Calabi-Yau (and other special holonomy) metrics are local
maxima for the Yamabe invariant.
\end{theo}

We also explore a remarkable connection between the stability
problem for $M$ and positive mass theorem on ${\mathbb R}^3 \times
M$ (One is reminded of Schoen's celebrated proof of the Yamabe
problem which makes essential use of the positive mass theorem).

This paper is organized as follows. We first review some
background on the variational problem. In fact, following
\cite{KW},  we use the first eigenvalue of the conformal Laplacian
instead of the total scalar curvature, which is essentially
equivalent but has the advantage of being conformally invariant in
a certain sense (Cf. Section 2). We then show how the parallel
spinor enables us to identify symmetric $2$-tensors with spinors
twisted by the cotangent bundle. This leads us to our first main
result, Theorem 1.1.

In Section 3 we prove the local stability theorem, Theorem 1.2.
This involves characterizing the kernel of Lichnerowicz Laplacian,
which is done according to the holonomy group. In each case, we
rely on the deformation theory of the special holonomy metric. For
example, in the Calabi-Yau case, we use the Bogomolov-Tian-Todorov
theorem \cite{t}, \cite{To} about the smoothness of the universal
deformation of Calabi-Yau manifold and Yau's celebrated soltuion
of Calabi conjecture \cite{y1,y2}. The deformation theory of $G_2$
metrics is discussed by Joyce \cite{J}. Essentially the same
techniques are employed in \cite{Wa2}.

We then discuss the connection with positive mass theorem which
presents a uniform approach. That is, one does not need separate
discussions for each of the special holonomy. The positive mass
theorem we use here is a special case of the result proved in
\cite{d}. We show that, if a metric of positive scalar curvature
is too close to the (Ricci flat) metric with nonzero parallel
spinor, then one can construct a metric on ${\mathbb R}^3 \times
M$ which has nonnegative scalar curvature, asymptotic to the
product metric at infinity, but has negative mass. Hence
contradictory to the positive mass theorem.

In the final section we make some remarks about compact Ricci flat
manifolds and point out the existence of scalar flat metrics which
are not Calabi-Yau on some Calabi-Yau manifolds. This existence
result depends on Theorem 1.4.

Our work is inspired and motivated by the recent work of
Hertog-Horowitz-Maeda \cite{HHM}.

{\em Acknowledgement:} The first and third authors are indebted to
Gary Horowitz for numerous discussions on his work \cite{HHM}. The
first author also thanks Thomas Hertog, John Lott, John Roe and
Gang Tian for useful discussion. The second author wishes to thank
Rick Schoen for stimulating discussions and encouragement.

\sect{The infinitesimal stability}

Let $(M, g_0)$ be a compact Riemannian manifold with zero scalar
curvature and $\op{Vol}(g_0)=1$. Our goal is to investigate the
sign of the scalar curvature of metrics near $g_0$. For any metric
$g$, we consider the conformal Laplacian $-\D_g+c_nS_g$, where
$c_n=(n-2)/4(n-1)$ and $S_g$ denotes the scalar curvature. Let
$\lam(g)$ be its first eigenvalue and $\ps_g$ the first
eigenfunction, normalized to satisfy $\int_M \ps_g dV_g=1$,  i.e.
\begin{gather}
-\D_g \ps_g+c_nS_g\ps_g=\lam(g)\ps_g,\label{eig} \\
\int_M \ps_gdV_g=1 \label{cons}
\end{gather}
So defined, $\ps_g$ is then uniquely determined and is in fact
positive.

The functional $\lam(g)$, studied first by Kazdan and Warner
\cite{KW}, has some nice properties. Though it is not conformally
invariant, its sign is conformally invariant. By the solution of
the Yamabe problem, any metric $g$ can be conformally deformed to
have constant scalar curvature. The sign of the constant is the
same as that of $\lam(g)$. In fact the metric $\ps_g^{4/(n-2)}g$
has scalar curvature $c_n\lam(g)\ps_g^{-4/(n-2)}$ whose sign is
determined by $\lam(g)$. Obviously $\lam(g_0)=0$ and the
corresponding eigenfunction is $\ps_0\equiv 1$. We present the
variational analysis of the functional $\lam$, essentially
following \cite{KW} with some modification and simplification.
(Our notations are also slightly different).

We first indicate that $\lam(g)$ and $\ps_g$ are smooth in $g$
near $g_0$. Let $U$ be the space of metrics on $M$, $V=\{\ps\in
C^{\infty}(M)|\int_M \ps dV_{g_0}=1\}$ and $W=\{\ps\in
C^{\infty}(M)|\int_M \ps dV_{g_0}=0\}$. We define $F:U\times V\ra
W$ as follows
\begin{equation}
F(g,\ps)=\left[-\D_g\ps^g+c_nS_g\ps^g-c_n\left(\int
S_g\ps^gdV_g\right)\ps^g\right]\frac{dV_g}{dV_{g_0}},
\end{equation}
where $\ps^g=\ps \frac{dV_{g_0}}{dV_g}$. Note $\int \ps^gdV_g=1$.
Obviously $F(g,\ps)=0$ iff $\ps^g$ is an eigenfunction of
$-\D_g+c_nS_g$ and is the first eigenfunction iff  $\ps>0$. We
have $F(g_0, \ps_0)=0$ and the linearization in the second
variable at $(g_0,\ps_0)$ is easily seen to be $-\D_{g_0}:W\ra W$.
Since this is an isomorphism, we conclude by the implicit function
theorem that $\lam(g)$ and $\ps_g$ are smoothly dependent on $g$
in a neighborhood of $g_0$.

Let $g(t)$ for $t\in (-\e,\e)$ be a smooth family of metrics with
$g(0)=g_0$. Before we analyze the variation of $\lam$ near $g_0$,
we collect a few formulas:
\begin{align}
\dot{\op{Ric}}&=\frac{1}{2}(\grd^*\grd h-2\cuv h)-\d^*\d h
-\frac{1}{2}D^2\op{tr}h+\op{Ric}\circ h, \label{dric} \\
\dot{S}&=-\lp h, \op{Ric}\rp+\d^2h-\D \op{tr}h, \label{dr}\\
\dot{\D}f&=-\lp h, D^2f\rp +\lp \d h+\frac{1}{2}d\op{tr}h, df\rp.
\label{dlap}
\end{align}
where $(\cuv h)_{ij}=R_{ikjl}h_{kl}$ denotes the action of the
curvature on symmetric $2$-tensors, $D^2$ denotes the Hessian, and
$k\circ h$ denotes the symmetric $2$-tensor associated to the
composition of $k$ and $h$ viewed as $(1,1)$-tensors via the
metric, i.e., as linear maps from $TM$ to itself. An ``upperdot''
denotes the derivative with respect to $t$ and $h=\dot{g}$.
Differentiating (\ref{eig}) in $t$ gives
\begin{equation}
\dot{\lam}\psi=-\D\dot{\psi}-\dot{\D}\psi+c_n(\dot{S}\psi+
S\dot{\psi})-\lam\dot{\psi}
\end{equation}
We integrate over $M$ and compute using the above formulas and
integration by parts
\begin{align*}
\dot{\lam}&=\int_M-\dot{\D}\ps+c_n(\dot{S}\ps+S\dot{\ps})-\lam\int_M
\dot{\ps}
\\
&=\int_M \lp h, D^2\ps\rp-\lp \d h+\frac{1}{2}d\op{tr}h,d\ps\rp +
c_n[-\lp h, \op{Ric}\rp \ps+\d^2h\ps-\D
\op{tr}h\ps+S\dot{\ps}]+\frac{\lam}{2}\int_M\ps\op{tr}h \\
&=c_n\int_M\lp
h,-\ps\op{Ric}+D^2\ps\rp+\frac{n}{n-2}\op{tr}h\D\ps+S\dot{\ps}
+\frac{\lam}{2}\int_M\ps\op{tr}h
\end{align*}
Therefore the first variation formula is
\begin{equation}\label{1v}
\dot{\lam}=c_n\int_M\lp h, -\ps\op{Ric}+D^2\ps+\frac{n}{n-2}\D\ps
g\rp+S\dot{\ps} +\frac{\lam}{2}\int_M\ps\op{tr}h.
\end{equation}
As $g_0$ is scalar flat, we have $\lam(0)=0, \ps_0=1$ and  hence
the elegant
\begin{equation}\label{1v0}
\dot{\lam}(0)=-c_n\int_M\lp\op{Ric}(g_0),h\rp dV_{g_0}.
\end{equation}
This shows that $g_0$ is a critical point of $\lam$ iff it is
Ricci flat.

As a corollary we have
\begin{prop}[Bourguignon] If $g_0$ has zero scalar curvature but non-zero Ricci
curvature, then it can
be deformed to a metric of positive scalar curvature.
\end{prop}

\Pf Take $h=-\op{Ric}(g_0)$ and $g(t)=g_0+th$. Then
$\dot{\lam}(0)=c_n\int_M|\op{Ric}(g_0)|^2dV_{g_0}>0$ and hence
$\lam(g(t))>0$ for $t>0$ small. Then $\ps_t^{4/(n-2)}g(t)$ has
positive scalar curvature with $\ps_t$ being the positive first
eigenfunction of $g(t)$. \qed

Now, a natural question is then, what happens if $\op{Ric}(g_0)$
is identically zero. This is exactly the question discussed by
Hertog, Horowitz and Maeda \cite{HHM} who, based on some physical
arguments, argued that in a neighborhood of a Calabi-Yau metric
there is no metric of positive scalar curvature. For this purpose
we need to derive the second variation for $\lam$.

We now assume that $g_0$ is Ricci flat. Differentiating (\ref{1v})
at $t=0$ and using (\ref{dric}) and (\ref{dr}) we get
\begin{align*}
\ddot{\lam}(0)&=c_n\int_M\lp h,
-\dot{\op{Ric}}+D^2\dot{\ps}+\frac{n}{n-2}\D\dot{\ps} g_0\rp
+\dot{S}\dot{\ps} \\
&=c_n\int_M\lp h, -\frac{1}{2}(\grd^*\grd h-2\cuv h)+\d^*\d
h+\frac{1}{2}D^2\op{tr}h+2D^2\dot{\ps} +\frac{2}{n-2}\D\dot{\ps}
g_0\rp  \\
&\stackrel{def}{=} c_n\int_M\lp h, \mathcal{F}(h)\rp.
\end{align*}
And $\dot{\ps}$ at $t=0$ appearing in the formula is determined by
the equation
\begin{equation}
\D\dot{\ps}=c_n\dot{S}=c_n(\d^2h-\D\op{tr}h).
\end{equation}
The symmetric tensor $h$ can be decomposed as $h=\bar{h}+h_1+h_2$,
where $h_1=L_Xg_0$ for some vector field $X$, $h_2=ug_0$ for some
smooth function $u$ and $\bar{h}$ is transverse traceless, that is
$\op{tr}\bar{h}=0$ and $\d\bar{h}=0$. Let $\ph_t$ be the flow
generated by the vector field $X$. Let $g(t)=\ph^*_tg_0$ and then
obviously its variation is $h_1=L_Xg_0$. Then
$\op{Ric}(g_t)=\ph^*_t\op{Ric}g_0=0$. Differentiating in $t$ and
using (\ref{dric}) we have
\begin{equation}
-\frac{1}{2}(\grd^*\grd h_1-2\cuv h)+\d^*\d
h_1+\frac{1}{2}D^2\op{tr}h_1=0.
\end{equation}
Similarly by working with the scalar curvature we have
\begin{equation}
\d^2h_1-\D\op{tr}h_1=0
\end{equation}
These two equation shows that $h_1$ has no contribution to
$\mathcal{F}(h)$. On the other hand it is easy to compute
\begin{gather*}
\op{tr}{h_2}=nu,\quad \d h_2=-du, \quad \d^*\d h_2=-D^2u \\
\grd^*\grd h_2-2\cuv h_2=-\D u g_0.
\end{gather*}
Therefore $\D\dot{\ps}=c_n(\d^2h_2-\D\op{tr}h_2)=-\frac{n-2}{4}\D
u$ and we may take $\dot{\ps}$ to be $-\frac{n-2}{4}u$. Putting
all these identities together we get
\begin{equation}
\mathcal{F}(h)=-\frac{1}{2}(\grd^*\grd \bar{h}-2\cuv \bar{h}).
\end{equation}
i.e. there is no contribution from $h_1$ and $h_2$. We summarize
our calculations as

\begin{prop}
Let $(M,g_0)$ be a compact Ricci flat manifold and $g(t)$ a smooth
family of metrics with $g(0)=g_0$ and $h=\frac{d}{dt}g(t)|_{t=0}$.
The second variation of $\lam$ at a Ricci flat metric is given by
\begin{equation}
\frac{d^2}{dt^2}\lam(g(t))|_{t=0} =-\frac{n-2}{8(n-1)}\int_M\lp
\grd^*\grd\bar{h}-2\cuv \bar{h},\bar{h}\rp dV_{g_0},
\end{equation}
where $\bar{h}$ the orthogonal projection of $h$ in the space of
transverse traceless symmetric $2$-tensors.
\end{prop}

\Rk The variational analysis for $\lam(g)$ parallels that for the
total scalar curvature functional as discussed in Schoen
\cite{Sch1}. The two functionals are essentially equivalent for
our purpose but $\lam(g)$ has the slight advantage that the second
variation on the space of conformal deformations is trivial.
\newline

The operator acting on symmetric 2-tensors
\begin{equation}
\mathcal{L}_gh=\grd^*\grd h-2\cuv h,
\end{equation}
appearing in the second variation formula for $\lam(g)$ is the so
called the Lichnerowicz Laplacian. To examine the nature of the
critical points of $\lam(g)$, it is important to determine if, for
Ricci flat metric, the Lichnerowicz Laplacian is nonnegative. This
was raised as an open question by Kazdan and Warner \cite{KW}.
Infinitesimally, the aforementioned work of Hertog, Horowitz and
Maeda indicates that it is expected to be so from physical point
of view at least for Calabi-Yau metrics.

The Lichnerowicz Laplacian is also of fundamental importance in
many other problems of Riemannian geometry. In general, however,
it is very difficult to study, for the curvature tensor is very
complicated. For metrics that are sufficiently pinched, there are
some results. See Besse \cite{Besse} 12.67. A very useful idea is
to view $h$ as $T^*M$-valued 1-form and then we have
\begin{equation}
(\delta^{\grd} d^{\grd}+d^{\grd}\delta^{\grd})h=\grd^*\grd h-\cuv
h+h\circ \op{Ric},
\end{equation}
where $d^{\grd}$ is the exterior differential operator on
$T^*M$-valued differential forms and $\delta^{\grd}$ its dual.
Therefore the left hand side is apparently positive semi-definite.
If $\op{Ric}=0$, the right hand side is different from
$\mathcal{L}_g$ only in the coefficients of the second term.
Though this formula does not help in general, it suggests that to
prove $\mathcal{L}_g$ to be positive semi-definite, one should
view $h$ as the section of some vector bundle with a differential
operator $P$ such that $\mathcal{L}_g=P^*P$. We show this is
possible for manifolds with parallel spinors.

We now assume $(M,g)$ is a compact spin manifold with the spinor
bundle $\Sr\ra M$. An excellent reference on spin geometry is
Lawson and Michelsohn \cite{LM}. Let $E\ra M$ be a vector bundle
with a connection. The curvature is defined as
\begin{equation}
R_{XY}=-\grd_X\grd_Y+\grd_Y\grd_X+\grd_{[X,Y]}.
\end{equation}
If $M$ is a Riemannian manifold, then for the Levi-Civita
connection on $TM$, we have $R(X,Y,Z,W)=\lp R_{XY}Z,W\rp$. We
often work with an orthonormal frame $\{e_1,\ldots,e_n\}$ and its
dual frame $\{e^1,\ldots,e^n\}$. Set
$R_{ijkl}=R(e_i,e_j,e_k,e_l)$.

The spinor bundle has a natural connection induced by the
Levi-Civita connection on $TM$. For a spinor $\s$, we have
\begin{equation}
R_{XY}\s=\frac{1}{4}R(X,Y,e_i,e_j)e_ie_j\cdot \s.
\end{equation}
If $\s_0\not=0$ is a parallel spinor, then
\begin{equation}\label{rxy=0}
R_{klij}e_ie_j\cdot \s_0=0.
\end{equation}
It is well known this implies $\op{Ric}=0$ by computation. From then on, we
assume $M$ has a parallel spinor $\s_0\not=0$, which, without loss
of generality, is normalized to be of unit length. We define a
linear map $\Phi: S^2(M)\ra \Sr\ts T^*M$ by
\begin{equation} \label{susy}
\Phi(h)=h_{ij}e_i\cdot \s_0\ts e^j.
\end{equation}
It is easy to check that the definition is independent of the
choice of the orthonormal frame $\{e_1,\ldots,e_n\}$.
\begin{lem}
The map $\Phi$ satisfies the following properties:
\begin{enumerate}
\item $\lp \Phi(h),\Phi(\tilde{h})\rp=\lp h,\tilde{h}\rp $, \item
$\grd_X\Phi(h)=\Phi(\grd_Xh)$.
\end{enumerate}
\end{lem}

\Pf We compute
\begin{align*}
\lp \Phi(h),\Phi(\tilde{h})\rp
&=h_{ij}\tilde{h}_{kl}\lp e_i\cdot\s_0\ts e^j,e_k\cdot\s_0\ts e^l\rp \\
&=h_{il}\tilde{h}_{kl}\lp e_i\cdot\s_0,e_k\cdot\s_0\rp \\
&=-h_{il}\tilde{h}_{kl}\lp \s_0,e_ie_k\cdot\s_0\rp \\
&=h_{kl}\tilde{h}_{kl} \\
&=\lp h,\tilde{h}\rp.
\end{align*}
This proves the first assertion. To prove the second one, we
choose our orthonormal frame such that $\grd e_i=0$ at $p$ and
compute at $p$
\begin{align*}
\grd_X\Phi(h)&=Xh_{ij}e_i\cdot \s_0\ts e^j \\
&=\grd_Xh(e_i,e_j)e_i\cdot \s_0\ts e^j \\
&=\Phi(\grd_X h).
\end{align*}
\qed

The following lemma is of crucial importance.
\begin{lem}
Let $h$ be a symmetric $2$-tensor on $M$. We have the following
Bochner type formula
\begin{equation}
\Dirac^*\Dirac\Phi(h)=\Phi(\grd^*\grd h-2\cuv h).
\end{equation}
\end{lem}

\Pf We choose an orthonormal frame $\{e_1,\ldots,e_n\}$ near a
point $p$ such that $\grd e_i=0$ at $p$. We compute at $p$
\begin{align*}
\Dirac^*\Dirac\Phi(h)&=e_k\cdot\grd_{e_k}(e_l\cdot\grd_{e_l}\Phi(h)) \\
&=\grd_{e_k}\grd_{e_l}h(e_i,e_j)e_ke_le_i\cdot\s_0\ts e^j \\
&=-\grd_{e_k}\grd_{e_k}h(e_i,e_j)e_i\cdot\s_0\ts e^j
  -\frac{1}{2}R_{e_ke_l}h(e_i,e_j)e_ke_le_i\cdot \s_0\ts e^j\\
&=\Phi(\grd^*\grd h)
   +\frac{1}{2}R_{kljp}h_{ip}e_ke_le_i\cdot\s_0\ts e^j
   +\frac{1}{2}R_{klip}h_{pj}e_ke_le_i\cdot\s_0\ts e^j. \\
\end{align*}
By using the Clifford algebra identity
$e_ie_j+e_je_i=-2\delta_{ij}$ twice we have
\begin{align*}
\frac{1}{2}R_{kljp}h_{ip}e_ke_le_i\cdot\s_0
&=-\frac{1}{2}R_{kljp}h_{ip}e_ke_ie_l\cdot\s_0-R_{kljp}h_{lp}e_k\cdot \s_0 \\
&=\frac{1}{2}R_{kljp}h_{ip}e_ie_ke_l\cdot\s_0+R_{kljp}h_{kp}e_l\cdot\s_0-R_{kljp}h_{lp}e_k\cdot
\s_0 \\
&=\frac{1}{2}R_{jpkl}h_{ip}e_ie_ke_l\cdot\s_0-2R_{kljp}h_{lp}e_k\cdot \s_0 \\
&=-2(\cuv h)_{kj}e_k\cdot\s_0,
\end{align*}
where in the last equality we used $R_{kljp}e_ke_l\cdot\s_0=0$ by
(\ref{rxy=0}). Similarly (in fact easier) one can show using also
the fact $\op{Ric}=0$
\begin{equation*}
\frac{1}{2}R_{klip}e_ke_le_i\cdot\s_0=0.
\end{equation*}
Thus we get
\begin{equation*}
\Dirac^*\Dirac\Phi(h)=\Phi(\grd^*\grd h-2\cuv h).
\end{equation*}
\qed

This yields
\begin{theo} \label{fmt}
The Lichnerowicz Laplacian $\mathcal{L}_g$ is positive
semi-definite and $\mathcal{L}_{g}h=0$ iff $\Dirac \Phi(h)=0$.
\end{theo}

Thus,
 \be W_{g}=\{h|\op{tr}_gh=0,\d h=0,\Dirac\Phi(h)=0\}. \ee
is the kernel of $\mathcal{L}_g$ on the space of transverse
traceless symmtric 2-tensors. We will study this finite
dimensional vector space in next section.

By working on a covering space, we also obtain
\begin{coro} If a compact Riemannian manifold $(M, g)$ has a cover which is spin
and admits nonzero parallel spinors, then the Lichnerowicz
Laplacian $\mathcal{L}_g$ is positive semi-definite.
\end{coro}

\Pf Let $\pi: \ (\Mh, \gh) \ra (M, g)$ be the cover. Clearly the
following diagram commutes
\begin{equation*} \begin{array}{cccc} \mathcal{L}_g:\ & S^2(M) & \ra & S^2(M) \\
& \pi^* \downarrow & & \pi^* \downarrow \\
\mathcal{L}_{\gh}:\ & S^2(\Mh) & \ra & S^2(\Mh). \end{array}
\end{equation*}

Now if we denote by $\lp \cdot \ , \ \cdot \rp$ the pointwise
inner product on symmetric $2$-tensors and $(\cdot \ , \ \cdot)$
the $L^2$ inner product, i.e., for example,
\[ ( h, h' )_g=\int_M \lp h, h' \rp_g dvol(g), \]
then we have

\be \lp \mathcal{L}_gh, h \rp_g=\lp \mathcal{L}_{\gh} \pi^*(h),
\pi^*(h) \rp_{\gh}. \ee

Thus, for a fundamental domain $F$ of $M$ in $\Mh$, one has

\be ( \mathcal{L}_gh, h )_g = \int_F \lp \mathcal{L}_{\gh}
\pi^*(h), \pi^*(h) \rp_{\gh} dvol(\gh). \ee

Since $(\Mh, \gh)$ has nonzero parallel spinor, we have $\Phi
(\mathcal{L}_{\gh}\pi^*h)=\Dirac^*\Dirac\Phi(\pi^*h)$ where the
map $\Phi$ is defined as in (\ref{susy}).

Without the loss of generality we take $\Mh$ to be the universal
cover. Since $M$ is Ricci flat, its fundamental group has
polynomial growth and therefore is amenable \cite{m}. Now we
choose $F$ as in \cite{b}. Namely, we pick a smooth triangulation
of $M$, and for each $n$-simplex in this triangulation, we pick
one simplex in $\Mh$ covering this simplex. We then let $F$ be the
union of all these simplices thus chosen in $\Mh$. Thus defined,
$F$ is a union of finitely many smooth $n$-simplices, but $F$ may
not be connected. With this choice of $F$, by Folner's theorem
\cite{b}, for every $\epsilon>0$, there is a finite subset $E$ of
the fundamental group such that the union of translates of $F$ by
elements of $E$,

\[ H= \bigcup_{g\in E} gF \]

satisfies

\[ \frac{{\rm area}(\partial H)}{{\rm vol}(H)} < \epsilon. \]

Hence

\ban ( \mathcal{L}_gh, h)_g & = & \frac{1}{\#E} \int_H
\lp \Phi(\mathcal{L}_{\gh} \pi^*(h)), \Phi(\pi^*(h))\rp_{\gh} dvol(\gh) \\
& = & \frac{\vol(M)}{\vol(H)} [ \int_H \lp \Dirac \Phi(\pi^*h),
\Dirac \Phi(\pi^*h)\rp_{\gh} dvol(\gh) + \int_{\partial H} \lp \nu
\cdot \Dirac \Phi(\pi^*h),
\Phi(\pi^*h)\rp_{\gh}\ int(\nu)dvol(\gh)]  \\
& \geq & - \frac{C \vol(M) {\rm area}(\pt H) }{\vol(H)}. \ean

Here we denote $\nu$ the outer unit normal of $\pt H$, and $C$
some constant depending on the $C^1$ norm of $h$ on $M$. Since the
right hand side of the last inequality above can be taken to be
arbitrarily small by appropriate choice of $E$, we obtain
\[ ( \mathcal{L}_gh, h)_g \geq 0. \]
\qed

\sect{The local stability theorem}

In this section we prove the following local stability theorem.

\begin{theo}\label{locfin}
Let $(M,g_0)$ be a compact, simply connected, irreducible
Riemannian spin manifold of dimension $n$ with a parallel spinor.
Then there exists a neighborhood $\mathcal{U}$ of $g_0$ in the
space of smooth Riemannian metrics on $M$ such that there exists
no metric of positive scalar curvature in $\mathcal{U}$.
\end{theo}

The key here is the identification of the kernel space
 \be W_{g}=\{h|\op{tr}_gh=0,\d h=0,\Dirac\Phi(h)=0\}. \ee
of $\mathcal{L}_g$ on the space of transverse traceless symmtric
2-tensors, according to the infinitesimal stability theorem. For
this purpose, we need to understand the geometry of $(M,g_0)$
better.

According to \cite{Wa}(cf. \cite{J} 3.6), if $(M,g_0)$ is a
compact, simply connected, irreducible Riemannian spin manifold of
dimension $n$ with a parallel spinor, then one of the following
holds
\begin{enumerate}
\item $(M, g_0)$ is flat,
\item $n=2m, m\geq 2$, the holonomy
group is $SU(m)$,
\item $n=4m, m\geq 2$, the holonomy group is
$Sp(m)$,
\item $n=8$, the holonomy group is $Spin(7)$,
\item
$n=7$, the holonomy group is $G_2$.
\end{enumerate}

In case 1, there is no metric of positive scalar curvature on $M$
by the work of Gromov-Lawson and Schoen-Yau (see \cite{LM}). In
cases 3 and 4, it is further shown in \cite{Wa} that the index of
the Dirac operator is nonzero, hence by Lichnerowicz's theorem
there is no metric of positve scalar curvature either. Therefore
the theorem is ``trivial'' execpt in cases 2 and 5.

Suppose $(M,g_0)$ is a compact Riemannian manifold of dimension
$n=2m$ with holonomy $SU(m)$. This is a Calabi-Yau manifold. By
Yau's solution of Calabi conjecture \cite{y1,y2} and the theorem
of Bogomolov-Tian-Todorov \cite{Bo}, \cite{t}, \cite{To}, the
moduli space $\Sigma$ of Calabi-Yau metrics is smooth of dimension
$h^{1,1}+2h^{m-1, 1}-1$ (it is one less than the usual number
because we normalize the volume and hence discount the trivial
deformation of scaling). Its tangent space at $g_0$ must be a
subspace of $W_{g_0}$ for $W_{g_0}$ is the Zariski tangent space
of the moduli space of Ricci flat metrics. In fact we have

\begin{lem} \label{iok}
$T_{g_0}\Sigma=W_{g_0}.$ \end{lem}

\Pf This follows from a theorem of Koiso \cite{K} which says
Einstein deformations of a K\"ahler-Einstein metric are also
K\"ahler, provided that first Chern class is nonpositive and the
complex deformation are unobstructed, which is guaranteed by
Bogomolov-Tian-Todorov theorem \cite{t}, \cite{To}. It can also be
easily seen form our approach. For a Calabi-Yau manifold its
spinor bunle is \be \Sr^+(M)=\bigoplus_{k\ even}
\wedge^{0,k}(M),\quad \Sr^-(M)=\bigoplus_{k\ odd} \wedge^{0,k}(M).
\ee The Clifford action at a point $p\in M$ is defined by \be
X\cdot\alpha
=\sqrt{2}(\pi^{0,1}(X^*)\wedge\alpha-\pi^{0,1}(X)\lrcorner \alpha)
\ee for any $X\in T_pM$ and $\alpha\in\Sr_p(M)$ And the parallel
spinor $\s_0\in C^\infty(\Sr^+(M))$ can be taken as the function
which is identically 1.

Let $J$ be the complex structure. Then we have $W_{g_0}=W^+\oplus
W^-$, where \be W^+=\{h\in W_{g_0}|h(J,J)=h\},\quad  W^-=\{h\in
W_{g_0}|h(J,J)=-h\} \ee We choose a local orthonormal $(1,0)$
frame $\{X_1,\ldots, X_m\}$ for $T^{1,0}M$ and its dual frame
$\{\theta^1,\ldots,\theta^m\}$. By straighforward computation we
have for $h\in W^+$ \be
\Phi(h)=h(\bar{X}_i,X_j)\bar{\theta}^i\ts\theta^j \ee which can be
identified with the real $(1,1)$ form
$\sqrt{-1}h(\bar{X}_i,X_j)\bar{\theta}^i\wedge\theta^j$. The Dirac operator
is then identified as
$\sqrt{2}(\overline{\partial}-\overline{\partial}^*)$ (cf. Morgan
\cite{mor}). Therefore $W^+$ is identified with the space of
harmonic $(1,1)$-forms orthogonal to the Kahler form $\omega$.
Similarly $W^-$ can be identified as $H^1(M,\Theta)-H^{0,2}(M)$,
where $\Theta$ is the holomorphic tangent bundle. As
$H^{0,2}(M)=0$ and $H^1(M,\Theta)\cong H^{m-1,1}(M)$ by the Hodge
theory, we have $\op{dim}W_{g_0}=h^{1,1}+2h^{m-1,1}-1$. This is
exactly the dimension of the moduli space of Calabi-Yau metrics.
\qed

We now turn to the proof of our local stability theorem in the
case of Calabi-Yau manifold.  Let $\mathcal{M}$ be the space of
Riemannian metric of volume $1$. By Ebin's slice theorem, there is
a real submanifold $\mathcal{S}$ containing $g_0$, which is a
slice for the action of the diffeomorphism group on $\mathcal{M}$.
The tangent space
\begin{equation}
T_{g_0}\mathcal{S}=\{h|\d_{g_0}h=0,
\int_M\op{tr}_{g_0}hdV_{g_0}=0.\}
\end{equation}
Let $\mathcal{C}\subset \mathcal{S}$ be the submanifold of constant scalar
curvatures metrics.If $g\in \mathcal{M}$ is a metric of positive
scalar curvature very close to $g_0$, then by the solution of the
Yamabe problem there is a metric $\tilde{g}\in \mathcal{C}$
conformal to $g$ and with constant positive scalar curvature.
Moreover as $g$ is close to $g_0$ which is the unique Yamabe
solution in its conformal class, $\tilde{g}$ is also close to
$g_0$. Therefore to prove the theorem, it suffices to work on
$\mathcal{C}$. It is easy to see \be
T_{g_0}\mathcal{C}=\{h|\d_{g_0}h=0, \op{tr}_{g_0}h=0.\} \ee

It contains the finite dimensional submanifold of Calabi-Yau
metrics $\mathcal{E}$. We now restric our function $\lam$ to
$\mathcal{C}$. It is identically zero on $\mathcal{E}$. Moreover by
Lemma \ref{iok}, $D^2\lam$ is negative definite on the normal
bundle. Therefore there is a possibly smaller neighborhood of
$\mathcal{E}\subset \mathcal{C}$, still denoted by $\mathcal{U}$, such
that $\lam$ is negative on $\mathcal{U}-\mathcal{E}$.

\bigskip
Next we consider the case 5, that is, the case of $G_2$ manifold.
Our basic references are Bryant \cite{br1,br2} and Joyce \cite{J}.
Let $(M,g_0)$ be a compact Riemannian manifold with holonomy group
$G_2$. We denote the fundamental $3$-form by $\ph$. With a local
$G_2$-frame $\{e_1,e_2,\ldots,e_7\}$ and the dual frame
$\{e^1,e^2,\ldots,e^7\}$ we have
\begin{align}
\ph&=e^{123}+e^{145}+e^{167}+e^{246}-e^{257}-e^{347}-e^{356}, \\
*\ph&=e^{4567}+e^{2367}+e^{2345}+e^{1357}-e^{1346}-e^{1256}-e^{1247}.
\end{align}
We also define the cross product $P:TM\times TM\ra TM$ by \be \lp
P(X,Y),Z\rp=\ph(X,Y,Z). \ee

The cross product has many wonderful properties. We list what we
need in the following lemma.

\begin{lem}\label{pp}
For any tangent vectors $X, Y, Z$
\begin{enumerate}
\item $P(X,Y)=-P(Y,X).$ \item $\lp P(X,Y),P(X,Z)\rp=|X|^2\lp
Y,Z\rp-\lp X,Y\rp \lp X, Z\rp, $ \item $P(X,P(X,Y))=-|X|^2Y+\lp
X,Y\rp X,$ \item $X\lrcorner(Y\lrcorner
*\ph)=-P(X,Y)\lrcorner\ph+X^*\wedge Y^*.$
\end{enumerate}
\end{lem}

\Pf The first three identities are proved in Bryant \cite{br1}.
The fourth can be proven by the same idea:  it is obviously true
for $X=e_1, Y=e_2$ and the general case follows by the
transitivity of $G_2$ on orthonormal pairs. \qed

The spinor bundle is $\Sr(M)=\mathbb{R}\oplus TM$ with the first
factor being the trivial line bundle. The Clifford action at $p\in
M$ is defined by \be X\cdot(a,Y)=(-\lp X,Y\rp, aX+P(X,Y)) \ee for
any $X,Y\in T_pM$. The parallel spinor $\s_0=(1,0)$. One easily
check $\ph(X,Y,Z)=-\lp X\cdot Y\cdot Z\cdot \sigma_0,\sigma_0\rp$.
It is also obvious that $\Sr(M)\ts T^*M=T^*M\oplus(TM\ts T^*M)$
and for any symmetric 2-tensor \be \Phi(h)=(0, h_{ij}e_i\ts e^j)
\ee We compute
\begin{equation}\label{dph}
\begin{split}
\Dirac{\Phi(h)}&=h_{ij,k}e_k\cdot(0,e_i)\ts e^j \\
&=h_{ij,k}(-\d_{ik},P(e_k,e_i))\ts e^j \\
&=(\d h,-h_{ij,k}P(e_i,e_k)\ts e^j).
\end{split}
\end{equation}

The $G_2$ structure gives rise to orthogonal decompostion of the
vector bundle of exterior differential forms. We are only
concerned with \be
\wedge^3(M)=\wedge^3_1(M)\oplus\wedge^3_7(M)\oplus
\wedge^3_{27}(M) \ee where
\begin{align}
\wedge^3_1&=\{a\ph \, |\ a\in \mathbb{R}\} \\
\wedge^3_7&=\{*(\ph\wedge\alpha)\, |\ \alpha\in T^*M\} \\
\wedge^3_{27}&=\{\alpha\in\wedge^3(M)\, | \
\alpha\wedge\ph=0,\alpha\wedge*\ph=0\}.
\end{align}
This also leads to the decomposition of the cohomology group \be
H^3(M,\mathbb{R})=\mathbb{R}\oplus H^3_7(M,\mathbb{R})\oplus
H^3_{27}(M,\mathbb{R}). \ee In fact one can show $
H^3_7(M,\mathbb{R})=0$. Therefore the Betti number
$b^3=1+b^3_{27}$, where $b^3_{27}=\op{dim}
H^3_{27}(M,\mathbb{R})$.

There is another natural isomorphism from the bundle of traceless
symmetric 2-tensors to $\wedge^3_{27}(M)$ \be \Psi:
S^2_0(M)\rightarrow \wedge^3_{27}(M) \ee defined by \be
\Psi(h)=h_{ij}e^i\wedge(e_j\lrcorner\ph). \ee It is proved by
Joyce (Theorem 10.4.4 in \cite{J}) that the moduli space of $G_2$
metrics is smooth and its tangent space at $g_0$ can be identified
with \be V_{g_0}=\{h\, |\, \op{tr}_{g_0}h=0,\d h=0,\Psi(h)\text{\
is harmonic}\}. \ee Moreover $\op{dim}V_{g_0}=b_3-1$. We have
$V_{g_0}\subset W_{g_0}$. To show they are equal, let $h\in
W_{g_0}$. By (\ref{dph}) we have \be\label{drh=0}
h_{ij,k}P(e_i,e_k)=0. \ee It suffices to show $\Psi(h)$ is
harmonic. We compute
\begin{align*}
d^*\Psi(h)&=-e_k\lrcorner \grd_{e_k}\Psi(h) \\
&=-h_{ij,k}e_k\lrcorner\left(e^i\wedge(e_j\lrcorner\ph)\right) \\
&=-h_{kj,k}e_j\lrcorner\ph+h_{ij,k}e^i\wedge\left(e_k\lrcorner(e_j\lrcorner\ph)\right) \\
&=(\d h)^{\sharp}\lrcorner\ph+h_{ij,k}e^i\wedge P(e_j,e_k)^* \\
&=0,
\end{align*}
where in the last step we used $\d h=0$ and (\ref{drh=0}).
Similarly
\begin{align*}
*d\Psi(h)&=*(e^k\wedge\grd_{e_k}\Psi(h)) \\
&=h_{ij,k}*\left(e^k\wedge e^i\wedge(e_j\lrcorner\ph)\right) \\
&=-h_{ij,k}e_k\lrcorner e_i\lrcorner(e^j\wedge *\ph) \\
&=-h_{ii,k}e_k\lrcorner*\ph+h_{ik,k}e_i\lrcorner*\ph
  -h_{ij,k}e^j\wedge\left(e_k\lrcorner (e_i\lrcorner *\ph)\right) \\
&=-h_{ij,k}e^j\wedge\left(e_k\lrcorner (e_i\lrcorner *\ph)\right),
\end{align*}
where in the last step we used $\op{tr}h=0$ and $\d h=0$. Then by
Lemma \ref{pp} we continue
\begin{align*}
*d\Psi(h)&=-h_{ij,k}e^j\wedge\left(e_k\lrcorner (e_i\lrcorner *\ph)\right) \\
&=-h_{ij,k}e^j\wedge \left(P(e_i,e_k)\lrcorner\ph\right)+h_{ij,k}e^j\wedge e^i\wedge e^k \\
&=0,
\end{align*}
where in the last step we used (\ref{drh=0}) and the fact that $h$
is symmetric. Therefore for any $h\in W_{g_0}$, the 3-form
$\Psi(h)$ is harmonic. This proves $W_{g_0}=V_{g_0}$ is the
tangent space to the moduli space of $G_2$ metrics. The rest of
argument is the same as in the Calabi-Yau case.

By the above argument we also obtain the following rigidity
theorem
\begin{theo} \label{rigidity}
Let $(M,g_0)$ be a compact, simply connected, irreducible
Riemannian spin manifold of dimension $n$ with a parallel spinor.
Then there exists a neighborhood $\mathcal{U}$ of $g_0$ in the
space of smooth Riemannian metrics on $M$ such that any scalar
flat metric in $\mathcal{U}$ must admit a parallel spinor (hence
Ricci flat in particular).
\end{theo}

\Pf We give the proof in the Calabi-Yau case. As shown in the
proof of Theorem \ref{locfin} $\lam$ is negative on
$\mathcal{U}-\mathcal{E}$, where $\mathcal{E}$ is the moduli space of
Calabi-Yau metrics. If $g\in \mathcal{U}$ is scalar flat then
$\lam(g)=0$. Therefore $g\in \mathcal{E}$.

\qed

Note that our proof actually gives a very nice picture of what
happens to scalar curvature near a metric with a parallel spinor
(let's call it special holonomy metric): one has the finite
dimensional smooth moduli of special holonomy metrics; along the
normal directions, the scalar curvature of the Yamabe metric in
the conformal class will go negative. That is, we have

\begin{theo} \label{locmax}
Let $(M,g_0)$ be a compact, simply connected, irreducible
Riemannian spin manifold of dimension $n$ with a parallel spinor.
Then $g_0$ is a local maximum of the Yamabe invariant.
\end{theo}

\Pf Recall that the Yamabe invariant of a metric $g$ is \be
\label{yi} \mu(g) = \inf_{f\in\, C^{\infty}(M),\ f>0} \frac{\int_M
(4\frac{n-1}{n-2} |df|_g^2 + S(g)f^2)dvol(g)}{(\int_M f^{2n/(n-2)}
dvol(g))^{(n-2)/n}}  \ee and it is a conformal invariant. The
corresponding Euler-Lagrange equation is \be \label{el}
-4\frac{n-1}{n-2} \D_gf + S(g)f=\mu(g) f^{(n+2)/(n-2)}. \ee Its
nontrivial solution, whose existence guaranteed by the solution of
Yamabe problem, gives rise to the so-called Yamabe metric
$f^{4/(n-2)} g$ which has constant scalar curvature $\mu(g)$. Note
that the left hand side of (\ref{el}) is (up to the positive
multiple of $4\frac{n-1}{n-2}$) the conformal Laplacian and the
numerator of the quotient in (\ref{yi}) is the quadratic form
defined by the conformal Laplacian (again up to the positive
multiple of $4\frac{n-1}{n-2}$). Thus, $\mu(g)<0$ if the first
eigenvalue of the conformal Laplacian $\lam(g)<0$. Now
$\mu(g_0)=\lam(g)=0$ since $g_0$ is scalar flat. Since we have
shown that $g_0$ is a local maximum of $\lam(g)$, our result
follows.  \qed

\sect{A uniform approach: connection with the positive mass
theorem}

In this section we explore a remarkable connection between the
stability of the Riemannian manifold $M$ and the positive mass
theorem on $\mathbb R^3 \times M$. This connection is pointed out
in the recent work \cite{HHM}. The positive mass theorem for
spaces which are asymptotic to $\mathbb R^k \times M$ at infinity
is proved in \cite{d}. The following result is a special case of
what is considered in \cite{d}.

\begin{theo} \label{pmt} Let $M$ be a compact Riemannian spin manifold with
nonzero parallel spinors. If $\tilde{g}$ is a complete Riemannian
metric on $\mathbb R^3 \times M$ which is asymptotic of order $>
1/2$ to the product metric at infinity and with nonnegative scalar
curvature, then its mass $m(\tilde{g}) \geq 0$. (Moreover,
$m(\tilde{g}) =0$ iff $\tilde{g}$ is isometric to the product
metric.)
\end{theo}

\Rk The result of \cite{d} is stated for simply connected compact
Riemannian spin manifolds $M$ with special holonomy with the
exception of quaternionic K\"ahler but the proof only uses the
existence of a nonzero parallel spinor, which is equivalent to the
holonomy condition by a result of \cite{Wa}. Also, the simply
connetedness of $M$ is not needed here since the total space is
the product $\mathbb R^3 \times M$.

Using the positive mass theorem, one can prove the following
deformation stability for compact Riemannian spin manifold with
nonzero parallel spinors, which is of course a special case of our
previous local stability theorem.

\begin{theo} \label{smt} Let $M$ be a compact spin manifold and $g$ a Riemannian
metric which admits nonzero parallel spinors. Then $g$ cannot be
deformed to metrics with positive scalar curvature. Namely, there
is no path of metrics $g_s$ such that $g_0=g$ and $S(g_s)>0$ for
$s>0$.
\end{theo}

\Rk The condition on the deformation can be relaxed to $S(g_s)
\geq 0$ and $S(g_{s_i})>0$ for a sequence $s_i \ra 0$.

This has the following interesting consequence, which, once again,
is a special case of the rigidity result, Theorem \ref{rigidity}.

\begin{theo} \label{sfdcy} Any scalar flat deformation of
Calabi-Yau metric on a compact manifold must also be Calabi-Yau.
\end{theo}

\Pf Let $(M, g)$ be a compact Calabi-Yau manifold and $g_s$ be a
scalar flat deformation of $g$, i.e., $S(g_s)=0$ for all $s$. We
will show by Theorem \ref{smt} that $g_s$ must also be Ricci flat.
It follows then by a theorem of Koiso and the
Bogomolov-Tian-Todorov theorem, $g_s$ must be in fact Calabi-Yau.

We first show that for $s$ sufficiently small, $g_s$ must be Ricci
flat. If not, there is a sequence $s_i \ra 0$ such that
$\Ric(g_{s_i}) \not= 0$. Hence this will also be true in a small
neighborhood of $s_i$. By Bourguignon's theorem, we can deform
$g_s$ slightly in such neighborhood so that it will now have
positive scalar curvature. Moreover we can do this while keeping
the metrics $g_s$ unchanged outside this neighborhood. This shows
that we have a path of metrics satisfying the conditions in the
Remark above. Hence contradictory to Theorem \ref{smt}.

Once we know that $g_s$ is Ricci flat for $s$ sufficiently small,
we then deduce that $g_s$ is Calabi-Yau for $s$ sufficiently
small, by using Koiso's Theorem \cite{Besse}, \cite{K} and
Bogomolov-Tian-Todorov Theorem \cite{t}, \cite{Bo}, \cite{To}.
Then the above argument applies again to extend to the whole
deformation. \qed

To prove Theorem \ref{smt}, we show that if there is such a
deformation of $g$, then one can construct a complete Riemannian
metric $\tilde{g}$ on $\mathbb R^3 \times M$ which is asymptotic
of order $1$ to the product metric at infinity and with
nonnegative scalar curvature, but with negative mass, hence
contradicting Theorem \ref{pmt}.

Now let $g(r)$, $r\ge 0$ be a one-parameter family of metrics on $M$. Consider
$\tilde{g}$ a warped product metric on $\mathbb R^3 \times M^k$ defined by
\be  \label{metric}
\tilde{g} = \left( 1-\frac{2m(r)}{r} \right)^{-1} dr^2 + r^2 ds^2 + g(r),
\ee
where $m(r)< r/2, m(0)=m'(0) =m''(0) =0$, $ds^2$ is the unit sphere metric on
$S^2$. The following lemma seems a folklore.

\begin{lem} If $g(r)$ is independent of $r$ for $r$ sufficiently large and
$\lim_{r \ra \infty} m(r) = m_{\infty}$ exists and finite, then $\tilde{g}$ on
$\mathbb R^3 \times M$ is asymptotically flat
of order $1$ with its mass $m(\tilde{g})=m_{\infty}$.
\end{lem}

We now turn to the scalar curvature of $\tilde{g}$. Given a point $(r,p,q) \in
\mathbb R^3 \times M^k$, choose a basis
 as follows:
$ U_0 = \left(1-\frac{2m(r)}{r}\right)^{\frac 12} \frac{\partial}{\partial r},
\ \{U_i\}$ an orthonormal basis of  $ds^2$ at $p$,
and $\{Y_\alpha\}$ a coordiante basis of $M$ at $q$.

Note that all brackets vanish except $[U_i,U_j]$ which belongs to $TS^2$, and
$\lp U_0,U_0\rp_{\tilde{g}} =1,\  \lp U_i,U_j\rp_{\tilde{g}} = r^2\delta_{ij}$.
Denote $\lp Y_\alpha,Y_\beta \rp_{\tilde{g}} = g_{\alpha \beta} (r)$.

By \cite[Formula (3.2)]{HHM}
\begin{lem} \label{scf}
The scalar curvature of $\tilde{g}$ at $(r,p,q)$ is given by
\be  \label{scalar}
\tilde{S}  = S_M (g(r)) + m' \left(\frac{4}{r^2} + \frac{g'}{r}\right) -
\frac{m}{r^2}g'  - \left(1-\frac{2m(r)}{r}\right) \left[ g'' + \frac{2}{r}g' +
\frac 14 (g')^2 -\frac 14 \partial_r g_{\alpha\beta} \partial_rg^{\alpha\beta}
\right],
\ee
where $g' = g^{\alpha\beta} \partial_r g_{\alpha\beta}, \ g'' = \partial_r g', \
S_M (g(r))$ is the scalar curvature of $(M, g(r))$ at $q$.
\end{lem}

This can be verified in local coordinates using Christoffel symbols. A simpler
way is to view
$S^2\times M^k$ with the product metric $r^2ds^2 + g(r)$ as a Riemannian
submanifold of $(\mathbb R^3 \times M^k, \tilde{g})$ and use Gauss equation to
compute some of the curvature tensors of
$\tilde{g}$. It seems this view point is useful whenever one has one parameter
family of metrics warping with $\mathbb R$. This type of warped metric appears
quite often, Cf (for example) \cite{W,BW},
\cite{P}. For this reason we present a proof here.

\Pf
Given tangent vectors $X,Y,Z,T$ at $(p,q) \in S^2 \times M^k$,
recall the Gauss equation \cite{doC}
\be  \label{gauss}
\lp \tilde{R}(X,Y)Z,T\rp = \lp R(X,Y)Z,T\rp - \lp B(Y,T), B(X,Z)\rp + \lp
B(X,T), B(Y,Z) \rp,
 \ee
 where $B$ is the second fundamental form.

 Note that
 \ban
 \lp B(U_i,U_j), U_0 \rp & = & \lp \tilde{\nabla}_{U_i}U_j, U_0\rp = -\frac12
U_0\lp U_i,U_j\rp =  -r\left(1-\frac{2m(r)}{r}\right)^{\frac 12} \delta_{ij}, \\
\lp B(Y_\alpha,Y_\beta), U_0 \rp & = &  -\frac12 U_0\lp Y_\alpha,Y_\beta\rp =
-\frac12
\left(1-\frac{2m(r)}{r}\right)^{\frac 12} g_{\alpha\beta,r}, \\
\lp B(U_i,Y_\beta), U_0 \rp & = & 0,
\ean
where $g_{\alpha\beta,r} = \frac{\partial g_{\alpha\beta}(r)}{\partial r}$.

Applying (\ref{gauss}) we get
\ba
\lp \tilde{R}(U_1,U_2)U_1,U_2\rp & = & r^2 - r^2 \left(1-\frac{2m(r)}{r}\right)
= 2mr  \label{uij} \\
\lp \tilde{R}(U_i,Y_\alpha)U_i,Y_\beta\rp & = & 0 -\frac12 r
\left(1-\frac{2m(r)}{r}\right)g_{\alpha\beta,r},  \label{uy}\\
\lp \tilde{R}(Y_\alpha,Y_\beta)Y_s,Y_l\rp & = & \lp
R(Y_\alpha,Y_\beta)Y_s,Y_l\rp_{(M,g(r))} - \frac 14
\left(1-\frac{2m(r)}{r}\right)\left[g_{\beta l,r} g_{\alpha s,r} -
g_{\alpha l,r}g_{\beta s,r}\right].  \label{yy} \ea

To compute the curvature tensors involving the $U_0$ direction,
note that \ban
\nabla_{U_0} U_0 & =& 0, \\
\nabla_{Y_\alpha} U_i & =& 0, \\
\nabla_{U_i} U_0 & =& \frac{1}{r}\left(1-\frac{2m(r)}{r}\right)^{1/2} U_i, \\
\nabla_{Y_\alpha} U_0 & =&\sum_{l} \left[ \sum_\beta \frac{1}{2}
\left(1-\frac{2m(r)}{r}\right)^{1/2} g_{\alpha\beta,r}g^{\beta l} \right] Y_l.
\ean
So
\ba
\lp \tilde{R}(U_0,U_i)U_0,U_i\rp & = & \lp \nabla_{U_i}\nabla_{U_0}U_0
- \nabla_{U_0}\nabla_{U_i}U_0 + \nabla_{[U_0, U_i]}U_0, U_i\rp
= - \lp \nabla_{U_0}\nabla_{U_i}U_0, U_i \rp  \nonumber
\\
& = & -\left(1-\frac{2m(r)}{r}\right)^{\frac 12}
\frac{d}{dr} \left[\frac{1}{r}\left(1-\frac{2m(r)}{r}\right)^{\frac 12}\right]
r^2 -\frac{1}{r^2}\left(1-\frac{2m(r)}{r}\right) r^2
\nonumber \\
& = & -\frac{m}{r} + m'.  \label{u0}
\ea
Similarly
\ba
\lp \tilde{R}(U_0,Y_\alpha)U_0,Y_\beta\rp
& = & \frac12 \frac{m' r-m}{r^2}g_{\alpha\beta,r} - \frac12
\left(1-\frac{2m(r)}{r}\right)\sum_{s,l}\frac{d}{dr} \left[g_{\alpha s,r}g^{ls}
\right] g_{l\beta} \nonumber \\
& & -\frac 14 \left(1-\frac{2m(r)}{r}\right) \sum_{s,l} g_{\alpha s,r} g^{ls}
g_{l\beta,r} .  \label{y0}
\ea
Let $\tilde{R}_{iljs} = \lp \tilde{R}(X_i,X_l)X_j, X_s \rp$ and $\tilde{g}^{sl}$
the inverse matrix of $\tilde{g}_{sl}$.
Since the Ricci tensor $\tilde{R}_{ab}$ is the trace of the curvature tensor and
$U_0,\ U_1, U_2$ are orthogonal
to each other in $\tilde{g}$, we have
\be
\tilde{R}_{ab} = \sum_{s,l} \tilde{R}_{albs} \tilde{g}^{sl}  = \tilde{R}_{a0b0}
+ \sum_i
\tilde{R}_{aibi} \frac{1}{r^2} + \sum_{\alpha, \beta} \tilde{R}_{a\alpha b\beta}
g^{\alpha\beta}.
\ee
Therefore using (\ref{u0}) (\ref{y0})
\ba
\tilde{R}_{00} & = & 2\left( -\frac{m}{r^3} + \frac{m'}{r^2} \right)
+\frac12 \frac{m' r-m}{r^2} g'-\frac12\left(1-\frac{2m(r)}{r}\right)g''
\nonumber \\
& & - \frac 14 \left(1-\frac{2m(r)}{r}\right) \sum_{\alpha,\beta,s,l}g_{\alpha
s,r} g_{l\beta,r} g^{ls}g^{\alpha \beta}. \label{r00}
\ea
Similarly by (\ref{u0}) (\ref{uij}) (\ref{uy}) and (\ref{y0}) (\ref{uy})
(\ref{yy})
\ba
\tilde{R}_{ii}  & =  & -\frac{m}{r} +m' +\frac{2m}{r} -\frac{r}{2}
\left(1-\frac{2m(r)}{r}\right) g', \label{rii} \\
\tilde{R}_{\alpha\beta} & = & \frac12 \frac{m' r-m}{r^2}g_{\alpha\beta,r} -
\frac12 \left(1-\frac{2m(r)}{r}\right)\sum_{s,l}\frac{d}{dr} \left[g_{\alpha
s,r}g^{ls} \right] g_{l\beta}  \nonumber \\
& & -\frac 14 \left(1-\frac{2m(r)}{r}\right) \sum_{s,l}g_{\alpha\beta,r}g_{sl,r}
g^{ls} -\frac1r
\left(1-\frac{2m(r)}{r}\right)g_{\alpha\beta,r} + R_{\alpha\beta} (M,g(r))
\label{rab}
\ea

Finally since the scalar curvature $\tilde{S}$ is the trace of the Ricci tensor
and $U_0,\ U_1, U_2$ are orthogonal, using (\ref{r00}) (\ref{rii}) (\ref{rab})
we have
\ban
\tilde{S} & =&  \sum_{ab} \tilde{R}_{ab} g^{ab} = \tilde{R}_{00} + \sum_i
\frac{1}{r^2} \tilde{R}_{ii} + \sum_{\alpha\beta}
\tilde{R}_{\alpha\beta} g^{\alpha\beta} \\
& = & \frac{4m'}{r^2} -\frac{m' r-m}{r^2} g'-\left(1-\frac{2m(r)}{r}\right)g''
  - \frac{2}{r}\left(1-\frac{2m(r)}{r}\right) g' + S_M (g(r)) \\
& & - \frac 14 \left(1-\frac{2m(r)}{r}\right)
\sum_{\alpha,\beta,s,l} \left[ g_{\alpha\beta,r} g_{sl,r} -g_{\alpha s,r}
g_{l\beta,r} \right] g^{\alpha\beta}g^{sl} \\
& = & S_M (g(r)) + m' \left(\frac{4}{r^2} + \frac{g'}{r}\right) -
\frac{m}{r^2}g'  - \left(1-\frac{2m(r)}{r}\right) \left[ g'' + \frac{2}{r}g' +
\frac 14 (g')^2 -\frac 14 \partial_r g_{\alpha\beta} \partial_rg^{\alpha\beta}
\right]
\ean
\qed

We now turn to the construction of asymptotically flat metrics of the type
(\ref{metric}) with nonnegative scalar curvature
and negative mass.

\begin{prop} \label{locsta}
Let $g_s, \ s \in [0,1]$, be a one-parameter family of metrics on $M$ such that
\[  S_M(g_1) \geq a_0>0,  \ \ \ \ S_M(g_0) \geq 0 . \]
If in addition
\be \label{cond1} \max | \partial_s
(g_s)_{\alpha\beta}
\partial_sg_s^{\alpha\beta} | \leq \frac{1}{200}, \ \ \ \max |g_s'|
\leq \frac{1}{200},  \ \ \  \max |g_s''|\leq \frac{1}{200} \ee
(here we did not try to get the optimal constants), and
\[ S^{-}_M(g_s) \stackrel{def}{=} \max(0, -S_M(g_s)) \leq \frac{a_0}{10} \]
for all $s \in [0, 1]$,  then there exists $m(r)$ with $m(0)=m'(0) =m''(0) =0$,
$\lim_{r \ra \infty} m(r)=m_{\infty}< 0$ and
$g(r)$ with $g(r)=g_0$ for $r$ sufficiently large, such that the metric
$\tilde{g}$ in (\ref{metric}) has scalar curvature
$\tilde{S} \ge 0$.
\end{prop}

\Pf  The function $m(r)$ and the metrics $g(r)$ will be constructed separately
on the intervals $[0, r_1]$,
$[r_1, r_2]$, $[r_2, r_3]$ and $[r_3, \infty)$, with $r_1, r_2, r_3$ to be
chosen appropriately.

First of all, define
\be
g(r)=g_1\ \ \mbox{ for} \ \ 0 \leq r \leq r_2.
\ee
 For the function $m(r)$, we let $m(r )= - \frac{a_0}{12} r^3$ when $0
\leq r \leq r_1$ and just require that $m'(r) \geq - \frac{a_0}{4}
r^2$  for $r_1 \leq r \leq r_2$, $m(r_2)=m(r_1)$ and
$m'(r_2)=-m'(r_1) = \frac{a_0}{4} r_1^2$. We note that there are
many such choices for $m(r)$ on the interval $[r_1, r_2]$ and
further $r_2$ can be taken arbitrarily close to $r_1$. For
definiteness we put $r_2=r_1 +1$.

In the region $0\leq r \leq r_2$, since the metrics $g(r)=g_1$
does not change with $r$, one easily sees from the scalar
curvature formula, Lemma \ref{scf}, that \be \tilde{S} = S_M(g_1)
+ m'(r) \frac{4}{r^2} \geq 0.           \ee

For the interval $[r_2, r_3]$, we define $m(r)$ to be the linear function
\[
m(r)= m(r_2) + \frac{a_0}{4} r_1^2 (r-r_2) \] and
\[ g(r)= g_{\frac{r_3-r}{r_3-r_2}}. \]
To make sure that the scalar curvature $\tilde{S}$ of the total
space remains nonnegative in this region, we need to choose the
parameters $r_1, r_3$ appropriately. Let $C_1= \max | \partial_s
(g_s)_{\alpha\beta} \partial_sg_s^{\alpha\beta} |$, $C_2= \max
|g_s'|$, $C_3=\max |g_s''|$. Then we have \ban \tilde{S} & \geq &
S_M(g(r)) +  \frac{a_0}{4} r_1^2 \left(\frac{4}{r^2} -
\frac{C_2}{r(r_3-r_2)}\right) - \frac{a_0}{12}
\frac{r_1^2}{r^2}\left| -r_1
+3(r-r_2) \right| \frac{ C_2}{r_3 -r_2} \\
& & - \left(1+ \frac{a_0}{6} \frac{r_1^2}{r} | -r_1 +3(r-r_2)
|\right) \left(\frac{C_1}{4(r_3-r_2)^2} + \frac{C_3}{(r_3-r_2)^2}
+ \frac{2C_2}{r(r_3-r_2)} + \frac{C_2^2}{4(r_3-r_2)^2}\right)  .
\ean This yields \be \label{bs} \tilde{S} \geq -S^{-}_M(g_s)+
\frac{a_0}{4} \frac{r_1^2}{r^2} (4 - A(r)) - \frac{1}{(r_3-r_2)^2}
B(r ), \ee where
\[ A(r)=C_2 \frac{r}{r_3-r_2} +\left( \frac{5C_2}{3} +\frac{4C_3 + C_1 + C_2^2}{6}
\frac{r}{r_3-r_2}\right) \frac{|-r_1 +3(r-r_2)|}{r_3-r_2} , \] and
\[ B(r)= \frac{C_1 + 4C_3 + C_2^2}{4} + 2C_2 \frac{r_3 -r_2}{r}.\]
We now set $r_3=r_2 + \frac{1}{7} r_1$. By the given bounds on
$C_1, C_2, C_3$ we have, for $r_2 \leq r \leq r_3, \ r_1 \ge 7$,
\[ |A(r)| \leq 3 , \ \ \ \ |B(r)| \leq 1\]
We then take $r_1$ sufficiently large so that
\[ \frac{1}{(r_3-r_2)^2} \leq \frac{a_0}{4} (\frac{1}{9})^2, \]
With these choices we find that $\tilde{S}\geq 0$ on the region $
r_2 \leq r \leq r_3$. Note that
\[ m(r_3)=m(r_2)+ \frac{a_0}{2} r_1^2(r_3-r_2)= -\frac{1}{84}r_1^3 a_0 \]
is still negative.

Finally, for $r \geq r_3$, we put $g(r )=g_0$ and choose $m(r)$ to
be increasing and limit to (say) $-\frac{1}{168}r_1^3 a_0$.
Clearly, $\tilde{S}\geq 0$ here. \qed

As a consequence, one has the following result.

\begin{coro} \label{infsta}
If $g_s, \ s \in [0,1]$ is a one-parameter family of metrics on $M$ such that
\[ S_M(g_0) \ge 0, \ \ \ \ S_M(g_s) > 0 \ \mbox{for}\ s \in (0,1), \]
then there exists $m(r)$ with $m(0)=m'(0) =m''(0) =0$, $\lim_{r \ra \infty}
m(r)=m_{\infty}< 0$ and $g(r)$ with $g(r)=g_0$
for $r$ sufficiently large, such that the
metric $\tilde{g}$ in (\ref{metric}) has scalar curvature $\tilde{S} \ge 0$.
\end{coro}

\Pf Let $\epsilon \in (0, 1)$ be a small constant. By taking $\epsilon$
sufficiently small, the hypothesis in  Proposition
\ref{locsta} can be obviously satisfied for the one-parameter family of metrics
$g_{\epsilon s}$ ($s\in [0, 1]$).
\qed

Theorem \ref{smt} follows from the Corollary above and the
Positive Mass Theorem.

In view of Proposition \ref{locsta} one also conclude that if $M$
is a compact spin manifold and $g$ a Riemannian metric which
admits nonzero parallel spinors, then there is no metrics of
positive scalar curvature in $U(g)$, where
\[ U(g)=\{ g'\in \mathcal{M}\ | \ \exists g_s \ {\mbox such \ that} \
g_0=g, g_1=g' \ {\mbox and} \ (\ref{cond1}) \ {\mbox holds};
S^-(g_s) \leq \frac{\min S(g')}{10} \  if \ S(g')>0. \}
\]

Condition (\ref{cond1}) essentially defines a scale invariant
$C^2$ neighborhood of $g$. The restriction on the negative part of
the scalar curvature of the path above is more subtle. Although we
believe that $U(g)$ is a neighborhood of $g$, we have not been
able to prove it.

Finally we remark that, just as in the previous section, the
results here still hold if $M$ may not be spin but there is a spin
cover of $M$ which admits nonzero parallel spinor. To see this, we
note that by Cheeger-Gromoll splitting theorem \cite{cg} the
universal cover $\Mh$ of $M$ is the product of a Euclidean factor
with a simply connected compact manifold. It is not hard to see
that the positive mass theorem as used here extends to this
setting.

\sect{Remarks on compact Ricci flat manifolds}

By Cheeger-Gromoll's splitting theorem \cite{cg}, we only have to
look at simply connected compact Ricci flat manifolds. Also,
without loss of generality, we can assume that they are
irreducible. Since all Ricci flat symmetric spaces are flat, by
Berger's classification of holonomy group \cite{Besse}, the
possible holonomy groups for a simply connected compact Ricci flat
 $n$-manifold are: $SO(n)$; $U(m)$ with $n=2m$; $SU(m)$;
$Sp(k)\cdot Sp(1)$ with $n=4k$; $Sk(k)$; $G_2$; and $Spin(7)$.
However, a simply connected Ricci flat manifold with holonomy
group contained in $U(m)$ must in fact have holonomy group
$SU(m)$; similarly, a simply connected Ricci flat manifold with
holonomy group contained in $Sp(k)\cdot Sp(1)$ must in fact have
holonomy group $Sp(k)$ \cite{Besse}. Thus,

\begin{prop} A simply connected irreducible compact Ricci flat $n$-manifold either
has holonomy $SO(n)$ or admits a nonzero parallel spinor.
\end{prop}

It is an interesting open question of what, if any, restriction
the Ricci flat condition will impose on the holonomy. All know
examples of simply connected compact Ricci flat manifolds are of
special holonomy and hence admit nonzero parallel spinors.

From our local stability theorem, scalar flat metrics sufficiently
close to a Calabi-Yau metric are necessarily Calabi-Yau
themselves. One is thus led to the question of whether there are
any scalar flat metrics on a compact Calabi-Yau manifold which are
not Calabi-Yau.

\begin{prop} If $(M, g)$ is a compact simply connected $G_2$
manifold or Calabi-Yau $n$-fold (i.e. $\dim_{\mathbb{C}}=n$) with
$n=2k+1$, then there are scalar flat metrics on $M$ which are not
$G_2$ or Calabi-Yau, respectively.
\end{prop}

\Pf By the theorem of Stolz \cite{S}, $M$ carries a metric of
positive scalar curvature (because of the dimension in the case of
$G_2$ manifold, and $\hat{A}=0$ in the case of Calabi-Yau
\cite{Wa}, \cite{J}). Let $g'$ be such a metric. Then, its Yamabe
invariant $\mu(g')>0$. On the other hand, every compact manifold
carries a metric of constant negative scalar curvature
\cite{Besse}, say $g''$. Then $\mu(g'') <0$. Now consider
$g_t=tg'+(1-t)g''$. Since the Yamabe invariant depends
continuously on the metric \cite{Besse}, there is a $t_0 \in
(0,1)$ such that $\mu(g_{t_0})=0$ and $\mu(g_t) >0$ for $0\leq t <
t_0$. The Yamabe metric in the conformal class of $g_{t_0}$ is
scalar flat but cannot be a $G_2$ (Calabi-Yau, resp.) by Theorem
\ref{locmax}. \qed

\Rk For the other special holonomy metric (including the
Calabi-Yau $SU(2k)$), the index $\hat{A}\not=0$ \cite{Wa},
\cite{J}. Thus there is no positive scalar curvature metric on
such manifold. Scalar flat metrics on such manifolds are
classified by Futaki \cite{fu}: they are all products of the
special holonomy metrics.

One is also led naturally to the following question: are there any
scalar flat but not Ricci flat metrics on a compact Calabi-Yau (or
$G_2$) manifold? A more interesting and closely related question
is:
\newline

\textbf{Question}: Is any Ricci flat metric on a compact
Calabi-Yau manifold necessarily K\"ahler (and hence Calabi-Yau)?
\newline

We note that the deformation analogue of this question has a
positive answer by the work of Koiso, Bogomolov, Tian, Todorov.
Moreover, in (real) dimension four, the answer is positive by the
work of Hitchin \cite{H2} on the rigidity case of Hitchin-Thorpe
inequality.

\end{document}